# Isaac Newton as a Probabilist

**Stephen M. Stigler**

*Abstract.* In 1693, Isaac Newton answered a query from Samuel Pepys about a problem involving dice. Newton's analysis is discussed and attention is drawn to an error he made.

On November 22, 1693, Samuel Pepys wrote a letter to Isaac Newton posing a problem in probability. Newton responded with three letters, first answering the question briefly, and then offering more information as Pepys pressed for clarification. Pepys (1633–1703) is best known today for his posthumously published diary covering the intimate details of his life over the years 1660–1669, but Newton would not have been aware of that diary. He would instead have known of Pepys as a former Secretary of Admiralty Affairs who had served as President of the Royal Society of London from 1684 through November 30, 1686, the same period when Newton's great *Principia* was presented to the Royal Society and its preparation for the press begun. But Pepys' letter did not concern scientific matters. He sought advice on the wisdom of a gamble.

## 1. PEPYS' PROBLEM

The three letters Newton wrote to Pepys on this problem, on November 26 and December 16 and 23, 1693, are almost all we have bearing on Newton and probability. Some of the letters were published with other private correspondence in Pepys (1825, Vol. 2, pages 129–135; 1876–1879, Vol. 6, pages 177–181) and more completely in Pepys (1926, Vol. 1, pages 72–94). The letters were cited in a textbook by Chrystal (1889, page 563), where he gave Pepys' problem as an exercise, but they were little known until they were brought to a wide public attention when selections were reprinted with commentary independently by Dan Pedoe (1958, pages 43–48), Florence David (1959; 1962, pages 125–129) and Emil D. Schell (1960). These authors and several others, notably Chaundy and Bullard (1960), Mosteller (1965, pages 6, 33–35) and Gani (1982) have discussed the problem Pepys posed and Newton's solution. Others accorded it briefer notice, including Sheynin (1971), who dismissively relegated it to a footnote; Westfall (1980, pages 498–499), who gave unwarranted credence to the excuse Pepys opened his first letter with, that the problem had some connection to a state lottery; and Gjertsen (1986, pages 427–428). But none of these or any other writer seems to have noted that a major portion of Newton's solution is wrong. The error casts an interesting light on how Newton thought about the matter, and it seems useful to revisit the question.

Since Pepys' original statement was, as Newton noticed, somewhat ambiguous, I will state the problem in paraphrase as it emerged in the correspondence:

Which of the following three propositions has the greatest chance of success?

A. Six fair dice are tossed independently and at least one "6" appears.

B. Twelve fair dice are tossed independently and at least two "6"s appear.

C. Eighteen fair dice are tossed independently and at least three "6"s appear.

As it emerged in the correspondence, Pepys initially thought that the third of these (C) was the most probable, but when Newton convinced him after repeated questioning by Pepys that in fact A was the most probable, Pepys ended the correspondence and announced he would, using Mosteller's (1965, page 35) colorful later term, welsh on a bet he had made.

*Stephen M. Stigler is Ernest Dewitt Burton Distinguished Service Professor of Statistics, Department of Statistics, University of Chicago, Chicago, Illinois 60637, USA e-mail: stigler@galton.uchicago.edu.*

## 2. NEWTON'S SOLUTION

Newton stated the solution three times during the correspondence: first he gave a simple logical reason for concluding that A is the most probable, then he reported a detailed exact enumeration of the chances in each of the three cases, and finally he returned to the logical argument and gave it in more detail.

Newton's exact enumeration was elegant and flawless; it is equivalent to the solution as might be presented in an elementary class today. Newton worked from first principles assuming no knowledge of the binomial distribution; we can now express what he found by this calculation in terms of a random variable $X$ with a Binomial $(N, p)$ distribution as follows:

A. $P(X \geq 1) = 31031/46656 = 0.665$ when $N = 6$ and $p = 1/6$.

B. $P(X \geq 2) = 1346704211/2176782336 = 0.619$ when $N = 12$ and $p = 1/6$.

C. Here Newton simply stated that, "In the third case the value will be found still less."

In fact,

$$P(X \geq 3) = 60666401980916/101559956668416$$
$$= 0.597$$

when $N = 18$ and $p = 1/6$, as another of Pepys' correspondents (a Mr. George Tollet) found after much labor, while trying to duplicate Newton's results (Pepys, 1926, Vol. 1, pages 92–94).

Pepys had originally thought that C was the most probable; Newton's logical arguments and his careful enumeration of chances pointed in the contrary direction. But while the conclusion Newton reached is correct, only the enumeration stands up under scrutiny. To understand why, it will help to develop a heuristic understanding of why A is the most probable.

## 3. A HEURISTIC VIEW

Pepys' problem amounts to a comparison of three Binomial $(N, p)$ distributions with $p = 1/6$, namely those with $N = 6$, 12 and 18. He desired a ranking of $P(X \geq Np)$ for the three cases. Now, in all Binomial distributions where the mean $Np$ is an integer, $Np$ is *also* the median of the distribution (and indeed the mode as well). This is *always* true, surprisingly even in cases like those under study here, where the distributions are quite skewed and asymmetric. This is a byproduct of a proof that for any $N$ and any $p$, the difference between the mean and median of a binomial distribution is strictly less than $\ln(2) < 0.7$ (Hamza, 1995). So when the mean $Np$ is an integer the two must agree, and this implies in particular that in all these cases,

$$P(X \geq Np) \geq \tfrac{1}{2} \quad \text{and} \quad P(X \leq Np) \geq \tfrac{1}{2},$$

and so in each case $P(X \geq Np)$ exceeds $1/2$ by a fraction of the probability $P(X = Np)$. In fact, in the cases Pepys considered we have to a fair approximation $P(X \geq Np) \approx 1/2 + (0.4)P(X = Np)$. The ranking Newton calculated then reflects the fact that the size of the modal probability for a binomial distribution, $P(X = Np)$, decreases as $N$ increases and the distribution spreads out, $p$ being held constant. Indeed, as De Moivre would find by the 1730s, $P(X = Np)$ is well approximated by $1/\sqrt{(2\pi Np(1-p))} \approx 1.07/\sqrt{N}$ when $p = 1/6$. So in particular, the probabilities in A, B, C are about $1/2 + (0.4)(1.07)/\sqrt{N}$, an approximation that would give values 0.67, 0.62, 0.60, which agree with the exact values to two places. Chaundy and Bullard (1960) provide a cumbersome rigorous proof that this sequence is decreasing, in some generality.

Note that this approximation depends crucially upon the probabilities $P(X \geq 1)$, $P(X \geq 2)$ and $P(X \geq 3)$ of A, B, C being $P(X \geq Np)$ [i.e. $P(X \geq E(X))$] for the three respective distributions, and the result depends upon this as well. Franklin B. Evans observed this sensitivity already in 1961, finding, for example, that $P(X \geq 1 | N = 6, p = 1/4) = 0.8220 < P(X \geq 2 | N = 12, p = 1/4) = 0.8416$ (Evans, 1961). That is, the ordering of A and B that Newton found for fair dice can fail for weighted dice, and indeed will tend to fail when $p$ is sufficiently greater than $1/6$, even though they be tossed fairly and independently.

## 4. NEWTON'S LOGICAL ARGUMENT

In his first letter to Pepys on November 26, 1693, Newton had been content to give a short logical argument for why the chance of A must be the largest. He dissected the problem carefully, and made it clear that the proposition required that in each case *at least* the given number of "6"s should be thrown. Newton then restated the question and gave an apparently clear argument as to why the chance for A had to be the largest:

"What is the expectation or hope of A to throw every time one six at least with six dyes?

"What is the expectation or hope of B to throw every time two sixes at least with twelve dyes?

"What is the expectation or hope of C to throw every time three sixes at least with 18 dyes?

"And whether has not B and C as great an expectation or hope to hit every time what they throw for as A hath to hit his what he throws for?

"If the question be thus stated, it appears by an easy computation that the expectation of A is greater than that of B or C; that is, the task of A is the easiest. And the reason is because A has all the chances of sixes on his dyes for his expectation, but B and C have not all the chances on theirs. For when B throws a single six or C but one or two sixes, they miss of their expectations." (Pepys, 1926, Vol. 1, 75–76; Schell, 1960)

Newton's conclusion was of course correct but the argument is not. It is easy for us to see that it cannot work because the argument applies equally well for weighted dice, and as we now know, the conclusion fails if, for example, $p$ is $1/4$. Any correct argument must explicitly use the fact that 1, 2, 3 are the expectations for A, B, C, and Newton's does not. His enumeration did do so, but A would equally well have "all the chances of sixes on his dyes" even if the chance of a "6" is $1/4$. Newton's proof refers only to the sample space and makes no use of the probabilities of different outcomes other than that the dice are thrown independently, and so it must fail. But Newton does casually use the word "expectations"; might he not have had something deeper in mind? His subsequent correspondence confirms that he did not.

In his third letter of December 23, 1693, Newton returned to this argument and expanded slightly on it. He personified the choices by naming the player faced with bet A "Peter" and the player faced with bet B "James." He then considered a "throw" to be six dice tossed at once, so then Peter was to make (at least) one "6" in a throw, while James was to make (at least) two "6"s in two throws.

Newton then wrote, "As the wager is stated, Peter must win as often as he throws a six [i.e., makes at least one "6" among the six dice], but James may often throw a six and yet win nothing, because he can never win upon one six alone. If Peter flings a six (for instance) four times in eight throws, he must certainly win four times, but James upon equal luck may throw a six eight times in sixteen throws and yet win nothing. For as the question in the wager is stated, he wins not upon every single throw with a six as Peter doth, but only upon every two throws wherein he throws at least two sixes. And therefore if he flings but one six in the two first throws, and one in the two next, and but one in the two next, and so on to sixteen throws, he wins nothing at all, though he throws a six twice as often as Peter doth, and by consequence have equal luck with Peter upon the dyes." (Pepys, 1926, Vol. 1, page 89; Schell, 1960)

Here we can see more clearly how Newton was led astray: Even though in the first letter he had carefully pointed out that "throwing a six" must be read as "throwing at least one six," here he confused the two statements. His argument might work if "exactly one six" were understood, but then it would not correspond to the problem as he and Pepys had agreed it should be understood. Indeed, Peter will not necessarily register a gain with every "6": if he has two or more in the first "throw" of six dice, he wins the same as with just one. Newton reduced the problem to single "throws" where each throw is a Binomial ($N = 6, p = 1/6$), and he lost sight of the multiplicity of outcomes that could lead to a win. Many of Peter's wins (those with at least two "6"s, which occurs in about 40% of the wins) would be wins for James as well. And in some of James's wins (those with at least two "6"s in one-half of tosses and none in the other half, about 28% of James's wins) Peter would not have done so well on "equal luck" (he would have won but half the time). Evidently to make Newton's argument correct would take as much work as an enumeration!

## 5. CONCLUSION

Newton's logical argument failed, but modern probabilists should admire the spirit of the attempt. It was a simple appeal to dominance, a claim that all sequences of outcomes will favor Peter at least as often as they will favor James. It had to fail because the truth of the proposition depends upon the probability measure assigned to the sequences and the argument did not. But this was 1693, when probability was in its infancy.

Why has apparently no one commented upon this error before? There are several possible explanations, and no doubt each held for at least one reader. (1) The letters were read superficially, with no attempt to parse the somewhat archaic language of the logical proof, which after all points in the right direction. (2) The language was puzzling and unclear to the reader (and Newton was not available to ask), but it was accepted since he was, after all, Isaac Newton, and the calculation clearly showed he was sound on the important fundamentals. (3) The reader may even have seen that it was not a satisfactory argument, but drew back from accusing Newton of error, particularly since he got the numbers right.

In a sense the argument is more interesting because it is wrong. Newton was thinking like a great probabilist—attempting a "eureka" proof that made the issue clear in a flash. When successful, this is the highest form of mathematical art. That it failed is no embarrassment; a simple argument can be wonderful, but it can also create an illusion of understanding when the matter is, as here, deeper than it appears on the surface. If Newton fooled himself, he evidently took with him a succession of readers more than 250 years later. Yet even they should feel no embarrassment. As Augustus De Morgan once wrote, "Everyone makes errors in probabilities, at times, and big ones." (Graves, 1889, page 459)